\documentclass[12pt,leqno]{article}
\begin{document}
\newtheorem{theorem}{Theorem}[section]
\newtheorem{prop}{Proposition}[section]
\newtheorem{defin}{Definition}[section]
\newtheorem{rem}{Remark}[section]
\newtheorem{corol}{Corollary}[section]
\newtheorem{example}{Example}[section]
\title{HAMILTONIAN STRUCTURES ON FOLIATIONS}
\author{{\small by}\vspace{2mm}\\Izu Vaisman}
\date{}
\maketitle
{\def\thefootnote{*}\footnotetext[1]%
{{\it 2000 Mathematics Subject Classification}: 53D17.
Hamiltonian structure. Poisson structure.}}
\begin{center} \begin{minipage}{12cm}
A{\footnotesize BSTRACT. We discuss hamiltonian structures of the
Gelfand-Dorfman complex of projectable vector fields and
differential forms on a foliated manifold. Such a structure
defines a  Poisson structure on the algebra of foliated functions,
and embeds the given foliation into a larger, generalized
foliation with presymplectic leaves. In a so-called tame case, the
structure is induced by a Poisson structure of the manifold.
Cohomology spaces and classes relevant to geometric quantization
are also considered.}
\end{minipage}
\end{center}
\vspace{5mm}
\begin{center}\section{Preliminaries}\end{center} Let $\mathcal{S}$ be a
moving body with supplementary physical characteristics, expressed
by scalar parameters, which have no impact on the motion but
depend on the latter. For instance, the temperature of a rigid
body which moves with high friction.

The mathematical model of such a system will consist of a {\it
configuration space} which is an $s$-dimensional differentiable
manifold $N$ endowed with a $p$-dimensional foliation
$\mathcal{G}$ such that the supplementary parameters are the
coordinates along the leaves of $\mathcal{G}$, and the position
coordinates are constant along these leaves. Then, the {\it phase
space} of $\mathcal{S}$ will be the total space $M$ of the
annihilator bundle $\nu^*\mathcal{G}\subseteq T^*N$ of the tangent
bundle $T\mathcal{G}$, and $M$ is endowed with the natural lift
$\mathcal{F}$ of $ \mathcal{G}$, which is such that the leaves of
$ \mathcal{F}$ are covering spaces of the leaves of $\mathcal{G}$
(e.g., see \cite{Mol}).

Since the motion does not depend on the supplementary parameters,
the hamiltonian function $H$ of the system will be a
$\mathcal{F}$-{\it foliated function} on $M$ i.e., a function
which is constant along the leaves of $ \mathcal{F}$. On the other
hand, since we want the motion to determine the time evolution of
the supplementary parameters, we should be able to define the
hamiltonian vector field of $H$ as a foliated vector field on the
phase space of $\mathcal{S}$.

Therefore, $(M,\mathcal{F})$ should be endowed with a {\it
generalized hamiltonian structure} that prescribes foliated
hamiltonian vector fields to foliated functions. The aim of this
paper is to initiate the study of such hamiltonian structures.

The generalized hamiltonian structures we need may be defined
within the general Gelfand-Dorfman scheme of hamiltonian
structures on complexes over a Lie algebra \cite{{Dor},{GD}}. For
convenience, we refer to such complexes as  {\it Gelfand-Dorfman
complexes} \cite{V2}, and recall their definition below.
\begin{defin} \label{GD-complex}
{\rm A {\it Gelfand-Dorfman complex} consists of:\\
i) a real Lie algebra $(\chi,[\;,\;])$;\\ ii) a cochain complex of
real vector spaces
$$\mathcal{C}=\;(\bigoplus_{k=0}^\infty\Omega^k,\;
d:\Omega^k\rightarrow\Omega^{k+1},\;d^2=0);$$ iii) mappings $X
\mapsto i(X)\in L_{\bf{R}}(\Omega^k,\Omega^{k-1})$,
($\Omega^{-1}:=0$; $:=$ denotes a definition), defined for
all $X\in\chi$ and $k=0,1,2,\ldots$, such that\\
a) if $\alpha\in\Omega^1$ and $i(X)\alpha=0$ for all $X\in\chi$
then $\alpha=0$;\\ b) if $L_X:=di(X)+i(X)d$ then}
\begin{equation} \label{GD-eqn}
i(X)i(Y)+i(Y)i(X)=0,\;i([X,Y])= L_Xi(Y)-i(Y)L_X.
\end{equation}  \end{defin}

Usually, one says that $\mathcal{C}$ is a {\it complex over}
$\chi$, and the mapping $X \mapsto i(X)$ encountered in Definition
\ref{GD-complex} may be seen as a {\it representation of $\chi$
on} $ \mathcal{C} $. This mapping also defines a pairing
$$<\alpha,X>=<X,\alpha>:=i(X)\alpha,\hspace{5mm}X\in\chi,\alpha\in
\Omega^1,$$ and, in particular, one denotes $Xf:=<df,X>$,
$f\in\Omega^0,X\in\chi$.

A linear mapping $H\in L_{\bf R}(\Omega^1,\chi)$ is said to be
{\it skew symmetric} if
\begin{equation} \label{skewsym}
<\alpha,H\beta>=-<\beta,H\alpha>,\hspace{5mm} \forall\alpha,\beta
\in\Omega^1. \end{equation}

The hamiltonian structures of a Gelfand-Dorfman complex are
defined by generalizing the notion of a Poisson bivector (e.g.,
\cite{V1}). For this purpose, one notices that the formula
\cite{{Dor},{GD}}
\begin{equation} \label{Schouten}
[H,K](\alpha,\beta,\gamma):=\sum_{Cycl(\alpha,\beta,\gamma)}
\{<KL_{H\alpha}\beta,\gamma>+<HL_{K\alpha}\beta,\gamma>\},
\end{equation} where $H,K\in L_{\bf R}(\Omega^1,\chi)$ are skew symmetric and
$\alpha,\beta,\gamma\in\Omega^1$, may be seen as defining a
bracket $$[H,K]\in L_{alt,{\bf R}}((\Omega^1)^3,\Omega^0),$$ which
is a generalization of the Schouten-Nijenhuis bracket of bivector
fields on manifolds. We call the bracket (\ref{Schouten}) the {\it
Gelfand-Dorfman bracket}. Then, one defines
\begin{defin} \label{Hamstr} {\rm A skew-symmetric homomorphism $H\in L_{\bf
R}(\Omega^1,\chi)$ which satisfies the {\it Poisson condition}
$[H,H]=0$ is called a {\it hamiltonian structure} on the
Gelfand-Dorfman complex $(\chi, \mathcal{C})$.}
\end{defin}

For a hamiltonian structure one defines the following
generalizations of
classical notions:\\
i) $\forall f\in \Omega^0$, $X_f:=H(df)\in\chi$ is the {\it
hamiltonian vector} of $f$;\\ ii) $\forall f,g\in \Omega^0$,
$\{f,g\}:=X_fg$ is the {\it Poisson bracket}; this bracket is
skew-symmetric because of (\ref{skewsym}),and it satisfies the
Jacobi identity because (\ref{Schouten}) yields
$$[H,H](df,dg,dh)=2\sum_{Cycl(f,g,h)}\{\{f,g\},h\};$$
iii) $\forall\alpha,\beta\in\Omega^1$, one has a $\Omega^1$-{\it
bracket}
\begin{equation} \label{formbrack}
\{\alpha,\beta\}:=
L_{H\alpha}\beta-L_{H\beta}\alpha-d<H\alpha,\beta>,
\end{equation}
with the particular case \begin{equation}\label{difbrack}
\{df,dg\}=d\{f,g\}. \end{equation}

The $\Omega^1$-bracket (\ref{formbrack}) may be defined for any
skew-symmetric mapping $H\in L_{{\bf R}}(\Omega^1,\chi)$, and it
satisfies the following fundamental identities \cite{{KM},{V21}}
\begin{equation} \label{identity1} <\gamma,H\{\alpha,\beta\}>=
<\gamma,[H\alpha,H\beta]>+\frac{1}{2}[H,H](\alpha,\beta,\gamma),
\end{equation}
\begin{equation} \label{identity2}
\sum_{Cycl(\alpha,\beta,\gamma)}<\{\{\alpha,\beta\},\gamma\},X> =
[H,L_XH](\alpha,\beta,\gamma)\end{equation}
$$+\frac{1}{2}\sum_{Cycl(\alpha,\beta,\gamma)}[H,H](\alpha,\beta,d<\gamma,X>),$$
where $\alpha,\beta,\gamma\in\Omega^1$, $X\in\chi$, and
\begin{equation} \label{derivLieH}
L_XH(\alpha):=[X,H\alpha]-H(L_X\alpha). \end{equation}

In the Hamiltonian case $[H,H]=0$, it follows from
(\ref{identity2}) that the $\Omega^1$-bracket is a Lie algebra
bracket. Furthermore, under the supplementary {\it regularity
hypothesis}: if $\forall\alpha\in\Omega^1$ $<\alpha,X>=0$  then
$X=0$ ($X\in\chi$), $H$ is a homomorphism of Lie algebras i.e.,
\begin{equation}\label{Hhom} H\{\alpha,\beta\}=[H\alpha,H\beta].
\end{equation}
On the other hand, even without the regularity hypothesis,
(\ref{identity1}) shows that if we ask $H\in L_{\bf
R}(\Omega^1,\chi)$ to be skew symmetric and satisfy (\ref{Hhom}),
$H$ is a hamiltonian structure.
\begin{center} \section{Hamiltonian structures of foliations}
\end{center} With the motivation of Section 1 in mind, let us consider an
arbitrary $n$-dimensional differentiable manifold $M$ (in the
present paper ``everything" is of differentiability class
$C^{\infty}$ ) endowed with a $p$-dimensional foliation $
\mathcal{F} $. An object of $M$ that projects to the space of the
leaves of $ \mathcal{F}$ is called either {\it projectable} or
{\it foliated}. We refer the reader to \cite{Mol} for all the
notions of foliation theory which we are going to use.

The Lie algebra $\chi_{ \mathcal{F}}$ of the $
\mathcal{F}$-foliated vector fields and the complex of projectable
differential forms $\Omega_{ \mathcal{F}}=
\bigoplus_{k=1}^{q}\Omega^k_{ \mathcal{F}}$ ($q:=n-p$), with the
usual exterior differential and contraction operators $i(X)$,
$X\in\chi_{ \mathcal{F}}$,  define a Gelfand-Dorfman complex
associated with the pair $(M,\mathcal{F})$. One might consider
general hamiltonian structures on this complex, but, such a
structure may have a non-local character. We avoid non-locality by
\begin{defin} \label{Hamfol} {\rm A
{\it hamiltonian structure} on (or of) the foliation $
\mathcal{F}$ is a vector bundle morphism $h:\nu^*
\mathcal{F}\rightarrow TM$ ($\nu\mathcal{F}=TM/T\mathcal{F}$ is
the transversal bundle of $ \mathcal{F}$) such that the induced
map of cross sections $H:\Omega^1_{
\mathcal{F}}\rightarrow\chi(M)$ ($\chi(M)$ is the space of all the
tangent vector fields of $M$) is a hamiltonian structure of the
Gelfand-Dorfman complex of $(M,\mathcal{F})$.}
\end{defin}

In particular, Definition \ref{Hamfol} implies that the morphism
$h$ is skew symmetric (i.e., it satisfies (\ref{skewsym})
pointwisely), and that the values of the mapping $H$ are in
$\chi_{\mathcal F}$.
\begin{example} \label{example0} {\rm Any skew symmetric $h\in L_{\bf
R}(\nu^*\mathcal{F},T\mathcal{F})$ may be seen as a trivial
hamiltonian structure of the foliation $ \mathcal{F}$. Indeed,
formula (\ref{Schouten}) shows that $[H,H]=0$ if the values of $H$
are vector fields tangent to $ \mathcal{F}$} \end{example}
\begin{example} \label{example1}
{\rm Let $P$ be a Poisson bivector field on the foliated manifold
$(M,\mathcal{F})$, such that for any foliated function $f\in
\Omega^0_\mathcal{F}$ the hamiltonian vector field $X^P_f$ is a
foliated vector field. Then, $h:=\sharp_P|_{\nu^*\mathcal{F}}$
($\sharp_P:T^*M\rightarrow TM$,
$<\sharp_P\alpha,\beta>:=P(\alpha,\beta)$) defines a hamiltonian
structure of the foliation $ \mathcal{F}$.}
\end{example}
\begin{example} \label{example2}
{\rm A bivector field $P$ is called a {\it transversal Poisson
structure} of $ \mathcal{F}$ if the bracket
$$\{f,g\}:=P(df,dg)\hspace{5mm}(f,g\in C^{\infty}(M))$$
makes $\Omega^0_{\mathcal{F}}$ a Poisson algebra \cite{MV}. In
this case, again, $h:=\sharp_P|_{\nu^*\mathcal{F}}$ is a
hamiltonian structure of $ \mathcal{F}$. Moreover, for any
hamiltonian structure $h$ of $ \mathcal{F}$ and any choice of a
decomposition $TM=E\oplus T\mathcal{F}$, the bivector field $P$
defined by
$$\sharp_P|_{E^*\approx\nu^*\mathcal{F}}=h,\;\sharp_P|_{T^*\mathcal{F}}=0$$
is a transversal Poisson structure of $ \mathcal{F}$.}
\end{example}

We also show how to express hamiltonian structures of a foliation
$ \mathcal{F}$ by means of {\it adapted local coordinates}
$(x^a,y^u)$, where $a=1,\ldots,q;\,u=q+1,\ldots,n$, and
$x^a=const.$ are the local equations of $ \mathcal{F}$. In order
to get an expression by tensors, we fix a decomposition
$TM=E\oplus T\mathcal{F}$ where \cite{V0}
\begin{equation} \label{normalbase}
E=span\{X_a :=\frac{\partial}{\partial
x^a}-t_a^u\frac{\partial}{\partial
y^u}\},\;T\mathcal{F}=span\{\frac{\partial}{\partial y^u}\},
\end{equation} for some local coefficients $t^u_a$ and with the
Einstein summation convention. The local bases of $TM$ defined by
(\ref{normalbase}) have the dual co-bases
\begin{equation} \label{normalcobase}
dx^a,\;\theta^u:=dy^u+t^u_adx^a,
\end{equation} and $\nu^*\mathcal{F}=span\{dx^a\}$.

Then, a skew-symmetric morphism $h:\nu^*\mathcal{F}\rightarrow TM$
has local equations
\begin{equation}\label{localh} h(dx^a)=h^{ab}X_b+
k^{au}\frac{\partial}{\partial y^u}.\end{equation} The components
$h^{ab}$ define a global cross section $W$ of $\wedge^2E$,
therefore, a global cross section of $\wedge^2\nu\mathcal{F}$,
which is independent on the choice of $E$, and the components
$k^{au}$ define a global cross section of $E\otimes T\mathcal{F}$.
The following assertion is obvious
\begin{prop} \label{prop28} The morphism $h$ defined by
{\rm (\ref{localh})} is a hamiltonian structure of $ \mathcal{F}$
iff the cross section $W$ with local components $h^{ab}$ is
foliated and defines a structure of Poisson algebra on $\Omega^0_{
\mathcal{F}}$.
\end{prop}

The Poisson bracket defined by $W$ on $\Omega^0_{ \mathcal{F}}$ is
of the local type, and it has the following interpretation. Let
$U$ be an open neighborhood of $M$ such that the manifold $N$ of
the slices of $ \mathcal{F}$ in $U$ exists, and let
$p:U\rightarrow N$ be the natural projection (constant along the
slices of $ \mathcal{F}$ in $U$). Then
$h_N(p(x)):=p_*(x)\circ(h|_U(x))\circ p^*(p(x))$, $x\in U$, is the
morphism $\sharp_{P_N}$ of a well defined Poisson bivector field
$P_{N}$ on $N$, which defines the same local Poisson brackets as
$W$. ($h_N$ is well defined since the values of the mapping $H$
defined by $h$ are foliated vector fields.)

Furthermore, any Poisson algebra structure of local type on
$\Omega^0_{ \mathcal{F}}$ is defined by a family of foliated
hamiltonian structures on $\mathcal{F}$. Indeed, the required
structure is equivalent to a foliated section $W$ of
$\bigwedge^2\nu( \mathcal{F})$, which satisfies the Poisson
condition $[W,W]=0$. Choose a decomposition $TM=E\oplus
T\mathcal{F}$, and, $\forall\alpha\in\Omega^1_{ \mathcal{F}}$,
define $h(\alpha)$ to be the unique vector of $E$ with projection
$\sharp_W\alpha$ on $\nu\mathcal{F}$. Since by (\ref{Schouten})
$$[H,H](\alpha,\beta,\gamma)=[W,W](\alpha,\beta,\gamma)\hspace{5mm}
(\alpha,\beta,\gamma\in\Omega^1_{ \mathcal{F}}),$$  $h$ is a
hamiltonian structure of $ \mathcal{F}$, and $h$ induces $W$.

More exactly, if $h_0$ is one of the foliated hamiltonian
structures which define $W$, the whole family which defines $W$ is
$h_0+k$, where $k\in L_{\bf R}(\nu^*\mathcal{F},T\mathcal{F})$ is
skew symmetric. This holds since for any hamiltonian structure $h$
of $\mathcal{F}$ and any skew symmetric $k\in L_{\bf
R}(\nu^*\mathcal{F},T\mathcal{F})$, the corresponding morphisms
$H,K$ of global cross sections satisfy the relation $[H,K]=0$ (see
(\ref{Schouten})).
\begin{prop} \label{prop22} For any hamiltonian structure $h$ on a foliation
$ \mathcal{F}$, the generalized distribution $
\mathcal{H}:=T\mathcal{F}+\mathcal{H}_0$ ($ \mathcal{H}_0:={\rm
im}\,h$) is a projectable, completely integrable distribution, and
its leaves are presymplectic manifolds with kernel $T\mathcal{F}$.
Furthermore, $h({\rm ann}\,T\mathcal{H})=\mathcal{H}_0\cap
T\mathcal{F}$ ({\rm ann} denotes the {\it annihilator} of a vector
space or bundle).
\end{prop} \noindent{\bf Proof.} We continue to use the previous
notation. Let $x_0\in U\subseteq M$ where $U$ is a neighborhood
such that $ \mathcal{F}|_U$ has a $q$-dimensional, transversal
submanifold $N$. Since the projection $p$ is a submersion, if
$L_{p(x_0)}$ is the symplectic leaf of the Poisson structure $P_N$
through $p(x_{0})\in N$, $\tilde L_{x_{0}}:=p^{-1}(L_{x_{0}})$ is
an integral submanifold of $ \mathcal{H}$ through $x_0$. The
existence of these integral submanifolds shows the complete
integrability of $ \mathcal{H}$. Projectability follows from the
fact that $ \mathcal{H}$ is spanned by the projectable vector
fields $H(\alpha)$, $\alpha\in\Omega^1_\mathcal{F}$, and
$\mathcal{H}$ projects onto the symplectic distribution of $P_N$.
The lift of the symplectic form of $L_{p(x_0)}$ by $p^*$ yields
the required presymplectic form of the corresponding leaf of
$\mathcal{H}$. Finally, notice that $\alpha\in{\rm
ann}\,\mathcal{H}$ iff $\alpha=p^*(\lambda)$ for some $\lambda\in
{\rm ker}\,\sharp_{P_N}$, and then $p_*h(\alpha)=0$. This implies
$h({\rm ann}\,\mathcal{H})\subseteq\mathcal{H}_0\cap
T\mathcal{F}$. On the other hand, if $h(\alpha)\in T\mathcal{F}$,
we must have $\alpha=p^*(\lambda)$ where $\lambda\in {\rm
ker}\,\sharp_{P_N}$, and this justifies the converse inclusion.
(All these also follow immediately from the local equations
(\ref{localh}) of $h$.) Q.e.d.

The distribution $ \mathcal{H}$ will be called the {\it
characteristic distribution} of the hamiltonian structure $h$, and
its leaves constitute the {\it presymplectic foliation}. The
hamiltonian structure $h$ of the foliation $ \mathcal{F}$ on $M$
will be called {\it transitive} if the characteristic distribution
is $ \mathcal{H}=TM$. In this case, Proposition \ref{prop22} tells
us that $M$ is a presymplectic manifold with the kernel foliation
$ \mathcal{F}$, and that $TM=\mathcal{H}_0\oplus T\mathcal{F}$.
The latter equality also shows that the corresponding local
Poisson structures $P_N$ are the symplectic reduction of the
presymplectic form of $M$. Conversely, if $M$ is a presymplectic
manifold with the presymplectic $2$-form $\sigma$, and if $E$ is a
complementary distribution of the kernel foliation $ \mathcal{F}$
of $\sigma$, there exists a well defined, transitive, hamiltonian
structure $h$ of $ \mathcal{F}$ such that $\mathcal{H}_0=E$ and
the local Poisson structures $P_N$ are the symplectic reductions
of $\sigma$.
\begin{example} \label{example3}
{\rm Let $ \mathcal{H}$ be a coisotropic foliation  of dimension
$n+k$ $(k\leq n)$ of a symplectic manifold $M$ of dimension $2n$,
with the symplectic form $\omega$. It is well known that the
$\omega$-orthogonal distribution of $ \mathcal{H}$ is tangent to a
foliation $ \mathcal{F}$, and that, $\forall x\in M$, there exist
local coordinates $(x^a,x^u,y^i)$ around $x$ such that
$a=1,\ldots, p:=n-k$, $u=p+1,\ldots,n$, $i=1,\ldots,n$,
$x^a=const.$ are the local equations of $ \mathcal{H}$, and the
symplectic form has the canonical expression
\begin{equation} \label{cansympl} \omega=\sum_{a=1}^{p} dx^a\wedge
dy^a +\sum_{u=p+1}^{n}dx^u\wedge dy^u. \end{equation} (This result
is a Lie's theorem \cite{Lib}.) The local equations of the
foliation $ \mathcal{F}$ are
$x^a=const.,\,x^u=const.,\,y^u=const,$, and the computation of the
hamiltonian vector field $X_f^{\omega}$ of an $
\mathcal{F}$-foliated function (via (\ref{cansympl})) shows that
$X_f^{\omega}$ is an $ \mathcal{F}$-foliated vector field tangent
to the leaves of $ \mathcal{H}$. Therefore,
$h:=-\flat_{\omega}^{-1}|_{\nu^*\mathcal{F}}$ is a hamiltonian
structure of the foliation $ \mathcal{F}$ with the presymplectic
foliation $ \mathcal{H}$. Moreover, in this case we have
$T\mathcal{F}\subseteq\mathcal{H}_0$.}
\end{example}
\begin{example} \label{example4}
{\rm Example \ref{example3} can be generalized as follows. Let
$(M,\omega)$ be an almost symplectic manifold (i.e., we ask
$\omega$ to be non-degenerate but not necessarily closed), and let
$ \mathcal{H}$ be a coisotropic foliation such that the pullback
of $\omega$ to every leaf of $ \mathcal{H}$ is closed on the leaf.
Then formula (\ref{cansympl}) is to be replaced by
\begin{equation} \label{canpres} \omega=\sum_{a=1}^{p} dx^a\wedge
\varpi^a +\sum_{u=p+1}^{n}dx^u\wedge dy^u, \end{equation} where
$\varpi^a$ are linearly independent, local, $1$-forms which
contain only the differentials $dy^a$. Now, we obtain the
foliation $ \mathcal{F}$ and its hamiltonian structure $h$ in the
same way as in the symplectic case.} \end{example}

We finish this section by a remark about the chosen definition of
the notion of a hamiltonian structure on a foliation.

If we start with the physical motivation of Section 1, and do not
think of Gelfand-Dorfman complexes a priori, the natural
definition of a generalized hamiltonian structure (g.h.s.) that
suites the problem is that of an ${\bf R}$-linear morphism of
sheaves
\begin{equation} \label{newhamiton}
\Phi:\underline{\Omega^0_\mathcal{F}}
\longrightarrow\underline{\chi},\hspace{5mm}f\mapsto X_f,
\end{equation} (underlining means passing to germs of
the corresponding type of objects), such that the bracket defined
by
\begin{equation} \label{newbrack}
\{f,g\}=X_fg,\hspace{5mm}f,g\in\underline{\Omega^0_\mathcal{F}},
\end{equation} makes $\underline{\Omega^0_\mathcal{F}}$ a Poisson
algebra sheaf.

In particular, the action of a hamiltonian vector field $X_f$ on
foliated functions $g\in\Omega^0_\mathcal{F}$ depends only on the
first jet $j^1f$. This is not enough to ensure that the g.h.s. has
local type. A natural condition for the latter property is to ask
$X_f=0$ for all $f\in\Omega^0_\mathcal{F}$ such that $j^1_xf=0$ at
each point $x\in M$. If the g.h.s. structure $\Phi$ satisfies this
locality condition, $\Phi$ is completely defined by local vector
fields
\begin{equation} \label{newlocalh}
X_{x^a}=h^{ab}X_b+k^{au}\frac{\partial}{\partial y^u},
\end{equation} that satisfy the conditions of Proposition
\ref{prop28}.

Therefore, the generalized hamiltonian structures of local type
are exactly the hamiltonian structures of foliations which we
defined earlier.
\begin{center}\section{Tame hamiltonian structures}
\end{center}
The Gelfand-Dorfman complex of a foliation does not satisfy the
regularity hypothesis formulated at the end of Section 1. The
equality $<\alpha,X>=0$, $\forall\alpha\in\Omega^1_{
\mathcal{F}}$, only implies $X\in\Gamma T\mathcal{F}$ ($\Gamma$
denotes the space of global cross sections). Therefore,
(\ref{Hhom}), or the equivalent property
\begin{equation} \label{Hhomeq} X_{\{f,g\}}=[X_f,X_g],\hspace{5mm}
 \forall f,g\in\Omega^0_{\mathcal{F}},
\end{equation} obtained by taking $\alpha=df,\beta=dg$, $f,g\in\Omega^0_{
\mathcal{F}}$ in (\ref{Hhom}), may not hold, and we shall define
\begin{defin} \label{strong}
{\rm A skew symmetric morphism $h:\nu^* \mathcal{F} \rightarrow
TM$ which satisfies condition (\ref{Hhomeq}) is a {\it strong
hamiltonian structure} on $ \mathcal{F}$.}
\end{defin}
\begin{rem} \label{prop26} {\rm If $h$ is a strong hamiltonian
structure, the sheaf $\underline{\nu^*\mathcal{F}}$ has a natural
structure of a sheaf of twisted Lie algebras {\rm\cite{KT}} over
$(\bf R, \underline{\Omega^0_\mathcal{F}})$, with the action of
germs $\alpha\in\underline{\nu^*_\mathcal{F}}$ defined as the
action of $H(\alpha)$.}
\end{rem}

Formula (\ref{identity1}) shows that a strong hamiltonian
structure is hamiltonian. The hamiltonian structures indicated in
Examples \ref{example1} and \ref{example3} are strong but, this is
not necessarily true for Examples \ref{example2} and
\ref{example4}. If $h$ is a strong hamiltonian structure, the
generalized distribution $ \mathcal{H}_0= {\rm im}\,h$ is
involutive. Conversely, if $ \mathcal{H}_0$ is involutive and if $
\mathcal{H}_0\cap T\mathcal{F}=0$, $h$ is a strong hamiltonian
structure (use (\ref{identity1})). These facts suggest
\begin{defin} \label{tame} {\rm A hamiltonian structure $h$ of a
foliation $\mathcal{F}$ is {\it transversal} (to $\mathcal{F}$) if
there exists a differentiable complementary distribution $E$ of
$T\mathcal{F}$ ($E\oplus T\mathcal{F}=TM$) such that $
\mathcal{H}_0\subseteq E$. The distribution $E$ will be called an
{\it image extension} of $h$. (It is possible to have more than
one image extension.) A transversal hamiltonian structure of $
\mathcal{F}$ is a {\it tame} structure if all the brackets of
differentiable vector fields that belong to $ \mathcal{H}_0$ are
contained in an image extension $E$.(In the tame case, only such
image extensions will be used.)}
\end{defin}

A tame hamiltonian structure is strong (see (\ref{identity1})),
and a transversal, strong hamiltonian structure is tame. The
condition $ \mathcal{H}_0\cap T\mathcal{F}=0$, which is implicit
in the definition of transversality, is equivalent to $h({\rm
ann}\,\mathcal{H})=0$ and also to the fact that the rank of the
morphism $h$ is equal to the rank of the Poisson structures
induced by $h$ on the manifolds of local slices of $ \mathcal{F}$.
(See Proposition \ref{prop22} and formula (\ref{localh}). This
condition is not enough for transversality. Indeed, there always
exists a smallest regular distribution $\bar\mathcal{H}_0$ which
contains the generalized distribution $ \mathcal{H}_0$ but, we may
have $\bar\mathcal{H}_0\cap T\mathcal{F}\neq0$.
\begin{example} \label{example5} {\rm Let $TM=F\oplus F'$ be
a locally product structure on the manifold $M$, and $
\mathcal{F}$ the foliations tangent to $F$. Assume that one has a
Poisson algebra structure of the local type on $\Omega^0_{
\mathcal{F}}$. Then, the hamiltonian structure $h$ which induces
the former and has its hamiltonian vector field in $F'$ is tame.
Indeed, $F'$ is an image extension of $h$ of the kind required for
tame structures. Notice also that a transitive, tame, hamiltonian
structure must be of the locally product type shown in the
example.}
\end{example}
\begin{prop} \label{prop31} Let $h$ be a transversal
hamiltonian structure of the foliation $ \mathcal{F}$ with image
extension $E$. Then $h$ is tame with image extension $E$ iff the
Nijenhuis tensor $N_E$ of the projection $p_E:TM\rightarrow TM$ of
$TM=E\oplus T\mathcal{F}$ onto $E$ satisfies the condition
\begin{equation} \label{Nijcond}
N_E(h\alpha,h\beta)=0,\hspace{1cm}\forall\alpha,\beta\in\nu^*_x
\mathcal{F},\;\forall x\in M. \end{equation} \end{prop}
\noindent{\bf Proof.} Following the general definition of a
Nijenhuis tensor e.g., \cite{{KM},{V2}} and since $p_E^2=p_E$, for
$X,Y\in\Gamma TM$, one has
\begin{equation}\label{defNij} N_E(X,Y)=[p_EX,p_EY]-
p_E[p_EX,Y]-p_E[X,p_EY]+p_E[X,Y].
\end{equation} Consider the local equations (\ref{localh}) of
$h$ using an image extension $E$, which implies that $k^{au}=0$.
Then, $h$ is tame iff
$$H(dh^{ab})=[H(dx^a),H(dx^b)],$$ which is equivalent to
\begin{equation}\label{tamelocal}
h^{ac}h^{be}\tau^u_{ce}=0,\hspace{3mm} \tau^u_{ce}:=
\frac{\partial t^u_c}{\partial x^e}-\frac{\partial t^u_e}{\partial
x^c} +t^v_c\frac{\partial t^u_e}{\partial y^v} -
t^v_e\frac{\partial t^u_c}{\partial y^v}. \end{equation} The
invariant meaning of (\ref{tamelocal}) is exactly (\ref{Nijcond}).
Q.e.d.

In the case of a transversal hamiltonian structure $h$ on a
foliated manifold $(M,\mathcal{F})$ it is possible to extend the
hamiltonian formalism in a way similar to what was done for
presymplectic manifolds in \cite{V01}.

Let us recall that, if $(M,\mathcal{F})$ is a foliated manifold
and if $E$ is a complementary distribution of $T\mathcal{F}$, the
use of the local bases (\ref{normalbase}), (\ref{normalcobase})
yields a bigrading of tensor fields and differential forms, with
the convention that the first degree is the $E$-degree and the
second is the $T\mathcal{F}$-degree \cite{V0}. For instance, a
differential $k$-form is of bidegree $(s,t)$ if its local
expressions contain $s$ forms $dx^a$ and $t$ forms $\theta^u$
$(s+t=k)$. Then, one has a decomposition
\begin{equation}\label{decompd}
d=d'_{(1,0)}+d''_{(0,1)}+\partial_{(2,-1)}, \end{equation} and
$d^2=0$ is equivalent to \begin{equation} \label{dsquarrezero}
\begin{array}{c}
d''^2=0,\;\partial^2=0,\;d'^2+d''\partial+\partial d''=0,\vspace{2mm}\\
d'd''+d''d'=0,\;\partial d'+d'\partial=0. \end{array}
\end{equation}

Now, we return to the transversal hamiltonian structure $h$ of $
\mathcal{F}$, and fix an image extension $E$ of $h$. Then the
corresponding section mapping $H$ is well defined for any
differential form $\alpha\in\Omega^{(1,0)}(M)$ of bidegree
$(1,0)$, and $H\alpha\in\Gamma E$. For any differentiable function
$f\in C^{\infty}(M)$, we can define the  hamiltonian vector field
$X'_f\in \Gamma E$  by
\begin{equation} \label{hamiltgen} X'_f=H(d'f) \end{equation}
and $\forall f,g\in C^{\infty}(M)$ we get an {\it extended Poisson
bracket}
\begin{equation} \label{Poissongen} \{f,g\}':=X'_fg=
<Hd'f,dg>=<Hd'f,d'g>=-\{g,f\}'.\end{equation}

Furthermore, if $X\in\Gamma E$ and $\alpha\in\Omega^{(1,0)}(M)$,
(\ref{decompd}) leads to
\begin{equation}\label{decompL} L_X\alpha=L'_X\alpha+L''_X\alpha,
\end{equation} where
\begin{equation} \label{decompL2} L'_X=i(X)d'+d'i(X),\;
L''_X=i(X)d''+d''i(X). \end{equation}

Accordingly, it is possible to extend the Gelfand-Dorfman bracket
(\ref{Schouten}) to arbitrary $(1,0)$-forms $\alpha,\beta,\gamma$
by \begin{equation} \label{extendedSchouten}
[H,K]'(\alpha,\beta,\gamma):=\sum_{Cycl(\alpha,\beta,\gamma)}
\{<KL'_{H\alpha}\beta,\gamma>+<HL'_{K\alpha}\beta,\gamma>\},
\end{equation} where $H,K$ are defined by skew symmetric morphisms
$h,k:\nu^*\mathcal{F}\rightarrow E$. A straightforward computation
shows that the extended bracket is trilinear over $C^{\infty}(M)$,
and for a hamiltonian structure $h$ we have $[H,H]'(\alpha,\beta,$
$\gamma)=0$ for any $\alpha,\beta,\gamma\in\Omega^{(1,0)}(M)$.

In particular, using (\ref{Poissongen}), (\ref{decompL2}), one
gets
\begin{equation} \label{Jacobi'}
[H,H]'(d'f,d'g,d'k)=2\sum_{Cycl(f,g,k)}[\{\{f,g\}',k\}'
+d'^2f(X'_g,X'_k)]=0. \end{equation}
\begin{prop} \label{th31} If $h$ is a tame hamiltonian
structure on $(M,\mathcal{F})$ the Poisson bracket $\{\;,\;\}'$
defines a Poisson structure on the manifold $M$. \end{prop}
\noindent{\bf Proof.} For any foliation and any choice of a
complementary distribution $E$ one gets
\begin{equation} \label{squarred'}
d'^2f(X,Y)=<d''f,N_E(X,Y)>,\hspace{5mm}\forall f\in
C^{\infty}(M),\,\forall X,Y\in\Gamma E, \end{equation} where $N_E$
is the Nijenhuis tensor (\ref{defNij}). Indeed, if $X,Y\in\Gamma
E$, (\ref{defNij}) yields
\begin{equation} \label{NpeE}
N_E(X,Y)=p_{T\mathcal{F}}[X,Y],\end{equation} where
$p_{T\mathcal{F}}$ denotes the projection onto the second term of
the decomposition $TM=E\oplus T\mathcal{F}$. On the other hand,
$$d'^2f(X,Y)=d(d'f)(X,Y)=XYf-YXf-<d'f,[X,Y]>$$
$$=[X,Y]f-(p_E[X,Y])f=<df,p_{T\mathcal{F}}[X,Y]>=
<d''f,p_{T\mathcal{F}}[X,Y]>.$$ Thus, (\ref{squarred'}) is
justified, and the conclusion follows from the characterization
(\ref{Nijcond}) of the tame hamiltonian structures and formula
(\ref{Jacobi'}). Q.e.d.

Theorem \ref{th31} tells us that a tame hamiltonian structure $h$
is defined by a usual Poisson structure $P$ on the foliated
manifold $(M,\mathcal{F})$. The hamiltonian vector fields of
foliated functions with respect to $h$ coincide with those with
respect to $P$, $\sharp_P|_{E^*}=h$ and
$\sharp_P|_{T^*\mathcal{F}}=0$. Thus, the tame hamiltonian
structures are included in Example \ref{example1}. But, not all
the structures of Example \ref{example1} are tame.

Similarly, it is possible to extend the bracket (\ref{formbrack})
of foliated $1$-forms to any $\alpha,\beta\in\Omega^{(1,0)}(M)$ by
\begin{equation} \label{formbrackgen}
\{\alpha,\beta\}':=L'_{H\alpha}\beta-L'_{H\beta}\alpha-d'<H\alpha,\beta>.
\end{equation}

From (\ref{formbrackgen}), it follows that $\forall f,g\in
C^{\infty}(M)$ one has
\begin{equation} \label{flin}
\{f\alpha,g\beta\}'=fg\{\alpha,\beta\}'+f(H(\alpha)g)\beta -
g(H(\beta)f)\alpha.\end{equation} In particular, we see that the
bracket (\ref{formbrackgen}) is skew symmetric because it is such
for foliated $1$-forms, where it reduces to (\ref{formbrack}).

Let us also evaluate the bracket (\ref{formbrackgen}) on an
argument $X\in\Gamma E$. First we define $L'_XH\in L_{\bf
R}(\Omega^{(1,0)}(M),\Gamma E)$ by
\begin{equation}\label{L'H} L'_XH(\alpha):=
p_E[X,H(\alpha)]-H(L'_X\alpha). \end{equation} Taking the
derivative of (\ref{skewsym}) in direction $X$, and with the
decomposition (\ref{decompL}), we see that $L'_XH$ is skew
symmetric. Then, if the derivatives $L'$ of (\ref{formbrackgen})
are replaced by $L-L''$ one gets
\begin{equation} \label{brackX} \{\alpha,\beta\}'(X)=
H(\alpha)i(X)\beta-H(\beta)i(X)\alpha-<\alpha,L'_XH(\beta)>.
\end{equation}
In particular, if $\alpha=d'f$, $\beta=d'g$ (\ref{brackX}) yields
\begin{equation} \label{brackd'}
\{d'f,d'g\}'=d'\{f,g\}+L_{X'_g}d''f-L_{X'_f}d''g. \end{equation}
The result follows by an easy computation which takes into account
the fact that the space of $(1,0)$-forms is the annihilator of
$E$.

If $L_{X'_g}d''f=0$ $\forall g\in C^\infty(M)$, we will say that
$f\in C^\infty(M)$ is a {\it distinguished function} \cite{V01},
and we will denote by $\Omega^0_d$ the space of distinguished
functions. For instance, any foliated function is distinguished
but, not conversely. By separating the $(1,0)$-term and the
$(0,1)$-term in the definition of a distinguished function, we see
that $f\in\Omega^0_d$ iff: {\it a}) $d'f$ is a foliated $1$-form,
and {\it b}) $\mathcal{H}\subseteq{\rm ker}\,d'^2f$. Formula
(\ref{Jacobi'}) shows that the extended Poisson bracket of
distinguished functions satisfies the Jacobi identity, and {\it
a}) implies that $\{f,g\}'\in\Omega^0_\mathcal{F}$ $\forall f,g\in
\Omega^0_d$. Therefore, $\Omega^0_d$ is a Poisson algebra and
$\Omega^0_\mathcal{F}$ is an ideal of the former. Furthermore,
(\ref{brackd'}) implies \begin{equation} \label{d'-compat}
\{d'f,d'g\}'=d'\{f,g\}',\hspace{5mm}\forall f,g\in\Omega^0_{d},
\end{equation} and, if we take $f,g\in\Omega^0_d,
k\in C^\infty(M)$ in (\ref{Jacobi'}) and use (\ref{NpeE}), we get
\begin{equation} \label{Xespecial} X'_{\{f,g\}'}=p_E[X'_f,X'_g]
\hspace{5mm}f,g\in\Omega^0_{d}.
\end{equation}
\begin{prop}\label{th32} Let $h$ be a tame hamiltonian
structure of the foliation $ \mathcal{F}$, $E$ an image extension
of $h$, and $P$ the Poisson structure defined by the brackets
$\{\;,\;\}'$. Then, the triple $ (\nu^*\mathcal{F},\{\;,\;\}',h)$,
with the bracket {\rm (\ref{formbrackgen})}, is a Lie subalgebroid
of the cotangent Lie algebroid $(T^*M,\{\;,\;\}_P, \sharp_P)$.
\end{prop}
\noindent{\bf Proof.} The bracket $\{\;,\;\}_P$ is given by
(\ref{formbrack}) with $H$ replaced by $\sharp_P$, and, since
$\sharp_P|_{E^*}=h$, we have $\forall
\alpha,\beta\in\Omega^1_\mathcal{F}$
$$\{\alpha,\beta\}_P= \{\alpha,\beta\}=\{\alpha,\beta\}'.$$
Then, (\ref{flin}) implies
$$\{f\alpha,g\beta\}_P=\{f\alpha,g\beta\}',\hspace{5mm}\forall f,g\in
C^\infty(M),\forall\alpha,\beta\in\Omega^1_\mathcal{F}.$$ Q.e.d.

Now, let us notice that there exist an inclusion and a splitting
morphism of Lie algebroids
\begin{equation}\label{inclusion} \iota:\nu^*\mathcal{F}
\hookrightarrow T^*M,\;\pi=p_{E^*}:T^*M\rightarrow
\nu^*\mathcal{F}\hspace{3mm}(\pi\circ\iota=id),
\end{equation} where $p_{E^*}$ is the projection onto $E^*$ in the
decomposition $T^*M=E^*\oplus T^*\mathcal{F}$.
\begin{prop} \label{propnou1} Under the hypotheses of Proposition \ref{th32}, the
projection $\pi$ induces an injection $\pi^*$ of the de Rham
cohomology of the Lie subalgebroid $\nu^*\mathcal{F}$ into the
Lichnerowicz-Poisson cohomology of $(M,P)$. For any complex vector
bundle S over $M$, the Lichnerowicz-Poisson Chern classes
$c^{LP}_k(S)$ belong to the image of the injection $\pi^*$.
\end{prop} \noindent{\bf Proof.} For the definition of the de Rham
cohomology of Lie algebroids, see \cite{KM}; the
Lichnerowicz-Poisson cohomology is the de Rham cohomology of the
cotangent Lie algebroid $T^*M$ of the Poisson manifold $(M,P)$
(e.g., \cite{V1}). These definitions show the existence of
homomorphisms
$$ \begin{array}{c}H^*_{{\rm de \,Rham}}(M)
\stackrel{j^*_1}{\rightarrow}H^*_{LP}(M,P)\stackrel{\iota^*}
{\rightarrow}H^*(\nu^*\mathcal{F}),\vspace{1mm} \\H^*_{{\rm de
\,Rham}}(M) \stackrel{j^*_2}{\rightarrow}
H^*(\nu^*\mathcal{F})\stackrel{\pi^*}{\rightarrow}
H^*_{LP}(M,P),\end{array}$$ where the morphisms are naturally
induced by $j_1=\sharp_P,j_2=h,\iota,\pi$. (For instance, at the
level of cochains we define
$$(j_2^*\lambda)(\alpha_1,...,\alpha_k)=
\lambda(H\alpha_1,...,H\alpha_k),\;(\lambda\in\Omega^k(M),
\alpha_1,...,\alpha_k\in\Gamma E^*),$$ etc.) The following
relations are obvious: $\iota^*\circ j_1^*=j^*_2$, $\pi^*\circ
j^*_2=j^*_1$, $\iota^*\circ\pi^*=id$. The last one shows that
$\pi^*$ is injective; the others were mentioned for a later
utilization.

Now, we remind that the Lichnerowicz-Poisson Chern classes are the
$j_1^*$-image of the real Chern classes. Representatives of
$c^{LP}_k(S)$ are obtained by evaluating Chern-Weil polynomials on
the curvature of an arbitrary contravariant derivative
$^P\hspace{-1mm}D$ on $S$ (i.e., a connection of the Lie algebroid
$T^*M$ on $S$) like in the usual Chern-Weil theory \cite{V1}. In
particular, if $^h\hspace{-1mm}D$ is a connection of the Lie
algebroid $\nu^*\mathcal{F}$ on $S$ then
\begin{equation} \label{connection} ^P\hspace{-1mm}D_\alpha s=\,^h
\hspace{-1mm}D_{\pi\alpha} s
\end{equation} is a contravariant derivative on $S$, and, if $C$
denotes curvatures, one has \begin{equation} \label{curvature}
C_{^P\hspace{-1mm}D}=\pi^* C_{^h\hspace{-1mm} D},\end{equation}
where $\pi^*$ is used at the level of cochains. Now, the same
procedure of evaluating Chern-Weil polynomials on curvature
applied to $C_{^h D}$ yields Chern classes $c^h_k(S)\in
H^{2k}(\nu^*\mathcal{F})$, which are the $j_2^*$-images of the
real Chern classes. Furthermore, (\ref{curvature}) shows that
$c^{LP}_k(S)=\pi^*c^h_k(S)$. Q.e.d.
\begin{corol} \label{corolar}Let $h$ be a tame hamiltonian
structure and let $P$ be the bivector field of the Poisson
brackets $\{\;,\;\}'$. Then, there exists a prequantization bundle
of the $h$-Poisson bracket iff $\iota^*[P]\in j_2^*(H^2(M,{\bf
Z}))$.
\end{corol} \noindent{\bf Proof.} $[P]\in H^2_{LP}(M,P)$ is the
cohomology class defined by the cocycle $P$. We refer the reader
to \cite{V1} for the geometric quantization theory involved in the
corollary. Since $P$ defines the same Poisson brackets as $h$, the
existence of a prequantization bundle implies $[P]=j_1^*(\zeta)$
for some $\zeta\in H^2(M,Z)$, which implies
$\iota^*[P]=j_2^*(\zeta)$. Conversely, if this condition is
satisfied, and if (as a consequence of (\ref{connection})) we see
the Kostant-Souriau prequantization formula as
\begin{equation}\label{prequant} \hat f(s)=\,^h\hspace{-1mm}
D_{d'f}s+2\pi\sqrt{-1}fs,\hspace{5mm}s\in\Gamma K ,\end{equation}
where $K$ is the required prequantization bundle, the Dirac
quantization principle implies that $c_1^h(K)=\iota^*[P]$. Since
we assumed that $\iota^*[P]$ is an integral cohomology class, $K$
exists. Q.e.d.

Now, let us consider the case of a transversal hamiltonian
structure $h$ on $(M,\mathcal{F})$, and fix an image extension
$E$. In this case, we may still see the cross sections of
$\wedge^k E$ as a kind of generalized cochains with a coboundary
$\delta^{(k)}=\delta$ defined by
\begin{equation} \label{coboundary} (\delta
Q)(\alpha_0,...,\alpha_k)=\sum_{i=0}^k (-1)^i
H(\alpha)(Q(\alpha_0,...,\hat
\alpha_i,...,\alpha_k))\end{equation}
$$+\sum_{i<j=1}^k (-1)^{i+j}Q(\{\alpha_i,\alpha_j\}',\alpha_0,...,\hat
\alpha_i,...,\hat \alpha_j,..., \alpha_k),$$ where $Q\in\wedge^k
E$, $\alpha_i\in \Gamma E^*$ ($i=0,...,k$), and the hat denotes
the absence of the corresponding argument.

If we denote $\delta^2=\delta^{(k+1)}\circ\delta^{(k)}$, a
straightforward computation yields \begin{equation}
\label{squarredelta} (\delta^2Q)(\alpha_0,..., \alpha_{k+1})
=\sum_{i<j=1}^{k+1}(-1)^{i+j}\Delta_h(\alpha_i,\alpha_j)(Q(\alpha_0,...,\hat
\alpha_i,\end{equation}$$...,\hat \alpha_j,...,\alpha_{k+1}))
+\sum_{i<j<k=2}^{k+1}
(-1)^{i+j+k}Q(\sum_{Cycl(i,j,k)}\{\alpha_k,\{\alpha_i,\alpha_j\}'\}',$$
$$\alpha_0,...,\hat \alpha_i,...,\hat \alpha_j,..., \hat
\alpha_k,...,\alpha_{k+1}),$$ where
\begin{equation} \label{X(fifj)}
\Delta_h(\alpha_i,\alpha_j):=H(\{\alpha_i,\alpha_j\}')
-[H(\alpha_i),H(\alpha_j)].
\end{equation} Since $\delta^2\neq0$, we can only define
the {\it twisted cohomology spaces} (e.g., \cite{V7})
\begin{equation} \label{twistedcoh} H_{tw}^k(h):=\frac{{\rm
ker}\,\delta^{(k)}} {{\rm im}\,\delta^{(k-1)}\cap{\rm
ker}\,\delta^{(k)}}.\end{equation} For instance, by
straightforward computations one gets $$ H_{tw}^0(h)=\{f\in
C^\infty(M)\;/\;X'_f=0\},\;H_{tw}^1(h)=\frac{\{Q\in\Gamma
E\;/\;L'_QH=0\}}{\{X'_f\;/\;f\in C^\infty(M), L'_{X'_f}=0\}}.
$$
But, if we define $W'\in\Gamma\wedge^2E$ by $
W'(d'f,d'g)=\{f,g\}'$, we do not get a cocycle since $$(\delta
W')(d'f,d'g,d'k)=-2\sum_{{\rm Cycl}(i,j,k)}\{\{f,g\}',k\}',$$ and
the Jacobi identity may not hold.

Several interpretations of twisted cohomology as a usual
cohomology exist (e.g., \cite{V7}). For instance, the subspaces
$\tilde\mathcal{C}^k(h)={\rm
ker}(\delta^{(k+1)}\circ\delta^{(k)})$ with the coboundary
$\delta$ constitute a usual cochain complex $
\tilde\mathcal{C}(h)$, and $H_{tw}^k(h)$ are the usual cohomology
spaces of $ \tilde\mathcal{C}(h)$.

On the other hand, since the Poisson bracket $\{\;,\;\}'$ defines
a representation of the Lie-Poisson algebra $\Omega_d^0$ of
distinguished functions on the space $\Omega^0_\mathcal{F}$ of
foliated functions, we get corresponding cohomology spaces
$H_d^*(h):=H^*(\Omega^0_d,\Omega^0_\mathcal{F})$. Then, the
cochains $$c(f_1,...,f_k)=Q(d'f_1,...,d'f_k),\;Q\in\Gamma\wedge^k
E, \,f_1,\ldots,f_k\in\Omega^0_\mathcal{F},$$ with values in
$\Omega^0_\mathcal{F}$ and the coboundary (\ref{coboundary})
define the cochain complex of projectable cross sections of
$\wedge E$ with the Lichnerowicz-like coboundary (see \cite{V1})
$\delta Q=-p_{\Gamma\wedge^{k+1}E}[W,Q]$ ($p$ denotes the
projection), where $W$ defines the $h$-Poisson bracket of foliated
functions, and $[\;,\;]$ is the Schouten-Nijenhuis bracket. We may
say that the cohomology spaces, say $H^*_{LPb}(h)$, of this
complex are the {\it basic Lichnerowicz-Poisson cohomology spaces}
of $h$. The restriction of the cochain $W'$ to distinguished
functions is $W$, and we have a {\it fundamental class} $[W]\in
H^2_{LPb}(h)$.

Now, remember that a foliated manifold also has {\it basic de Rham
cohomology spaces} $H^*_b(M,\mathcal{F})$ \cite{Mol}, defined as
the cohomology spaces of the complex $(\Omega^*_\mathcal{F},d)$,
and there exist natural homomorphisms
$$\varphi:H^*_b(M,\mathcal{F})\rightarrow H^*_{{\rm
de\,Rham}}(M),\;\psi:H^*_b(M,\mathcal{F})\rightarrow
H^*_{LPb}(h),$$ induced by inclusion and $h$, respectively.

These facts have the following consequences for geometric
quantization. Assume that $[W]=\psi[\Phi]$ where $\varphi[\Phi]$
is an integral de Rham cohomology class. Then
$\Phi\in\Omega^2_\mathcal{F}$ is a closed $2$-form with integral
periods, such that \begin{equation} \label{integral}
\{f,g\}'=\Phi(X'_f,X'_g),\hspace{5mm}\forall f,g\in\Omega^0_d.
\end{equation} Accordingly, there exists a Hermitian line bundle
$K$ over $M$ with a connection $\nabla$ of curvature
$2\pi\sqrt{-1}\Phi$, and the Kostant-Souriau formula
\begin{equation} \label {KSultim} \hat fs=\nabla_{X_f}s+
2\pi\sqrt{-1}fs \end{equation} provides a prequantization such
that the Dirac principle holds for distinguished functions but,
generally, not for arbitrary functions (use (\ref{integral})). The
transitive case, i.e., presymplectic manifolds, was discussed in
\cite{V01}. {\small
\begin{center} \end{center}}\hspace*{7.5cm}{\small \begin{tabular}{l}
Department of Mathematics\\ University of Haifa, Israel\\ E-mail:
vaisman@math.haifa.ac.il \end{tabular}}
\end{document}